\documentclass[12pt,oneside]{article}

\usepackage{amscd}
\usepackage{amssymb,amsmath}
\usepackage[latin1]{inputenc}
\usepackage{verbatim}
\usepackage{multibib}

\def\Z{\mathbb{Z}}

\def\ZZ{\Z \oplus \Z}
\def\R{\mathbb{R}}
\def\C{\mathbb{C}}
\def\H{\mathbb{H}}

\def\P{\pi_{1}}

\def\G{\Gamma}

\newtheorem{thm}{Theorem}

\newtheorem{conj}{Conjecture}

\newtheorem{prop}{Proposition}[section]

\date{}

\title{A survey on Seifert fibre space conjecture}

\newcites{chrono}{References of the proof of the SFSC}
\newcites{gen}{References}

\begin{document}
\maketitle 
\begin{center}
{\sc Jean-Philippe PR\' EAUX}\footnote[1]{Centre de recherche de l'Arm\'ee de l'Air, Ecole de l'air, F-13661 Salon de
Provence air}\ \footnote[2]{Centre de Math\'ematiques et d'informatique, Universit\'e de Provence, 39 rue
F.Joliot-Curie, F-13453 marseille
cedex 13\\
\indent {\it E-mail :} \ preaux@cmi.univ-mrs.fr\\
{\it Mathematical subject classification : 57N10, 57M05, 55R65.\\
Preprint version 2010}}
\end{center}
\begin{abstract}
We recall the history of the proof of Seifert fibre space conjecture, as well as it motivations and its several
generalisations.
\end{abstract}
\addcontentsline{toc}{subsection}{Table of contents}
\tableofcontents

\subsection*{Introduction}\addcontentsline{toc}{subsection}{Introduction}
The reader is supposed to be familiar with general Topology, Topology of 3-manifolds, and algebraic Topology. Basic
courses can be found in celebrated books\footnote[1]{\ J.\textsc{Hempel}, \textsl{3-manifolds}, Princeton University Press, 1976.\\
\null\qquad H.\textsc{Seifert} and W.\textsc{Threlfall}, \textsl{A Textbook of Topology}, Academic Press, New-York, 1980.\\
\null\qquad E.\textsl{Spanier}, \textsl{Algebraic Topology}, Mc Graw-Hill, 1966.
}.

We work in the PL-category. We denote by a {\it 3-manifold} a PL-manifold of dimension 3 with
boundary (eventually
empty), which is moreover, unless otherwise stated, assumed to be connected. The boundary of a 3-manifold $M$ is denoted by $\partial M$.\\

In the topology of low-dimensional manifolds ({\it i.e.} of dimension at most 3) the fundamental group (or $\pi_1$) plays a central key-role. On the one hand several of the main topological properties of 2 and 3-manifolds can be rephrased in term of properties of the fundamental group and on the other hand in the generic cases the $\pi_1$ fully determines their homeomorphic type. That the $\pi_1$ generally determines their homotopy type follows from the fact that they are generically Eilenberg-McLane spaces (or $K(\pi,1)$, {\it i.e.} their universal cover is contractile); that the homotopy type determines generically the homeomorphism type appears as a rigidity property, or informally that the lack of dimension prevents the existence of too many 2 and 3-manifolds; this contrasts with higher dimensions.

It has been well known since a long time for 2-manifolds (or surfaces); advances in the study of 3-manifolds show that it seems to remain globally true in dimension 3.  For example, among several others,  one can think at the Poincar\'e conjecture, at the Dehn-loop and the sphere theorems of Papakyriakopoulos, at the torus theorem, at the rigidity theorem for Haken manifold, and at the Mostow's rigidity theorem for hyperbolic 3-manifolds, {\it etc...}. It provides (among other properties such as the equivalence of PL, topological and differential structures on 2, 3-manifolds) a somehow common paradigm for their study, much linked to combinatorial and geometric group theory, which has constituted into a particular discipline, {\sl low-dimensional topology}, among the the more general {\sl topology of manifolds}.\\

Among all 3-manifolds, a particular class introduced by H.Seifert in 1933, known as {\sl Seifert manifolds} or {\sl Seifert fibre spaces} has been since widely studied, well understood, and is having a great impact for understanding 3-manifolds. They suit many nice properties, whom majority were already known since the deep work of Seifert; nevertheless one of their major property, of major importance in the understanding of compact 3-manifolds, the so called {\sl Seifert fibre space conjecture}, has been a long standing conjecture before the proof was completed by the huge collective work involving several mathematicians,  among who F.Waldhausen, C.Gordon \& W.Heil, W.Jaco \& P.Shalen, P.Scott, G.Mess, P.Tukia, A.Casson \& D.Jungreis, and D.Gabai, for about half a century ; one of these 'monster theorems' which appeared in the twentieth century. It has turned to become another example
of the characteristic meaning of the $\pi_1$ for 3-manifolds.

The {Seifert fibre space
conjecture} characterizes those Seifert fibre spaces with infinite $\pi_1$ in the class of oriented irreducible 3-manifolds in term of a property
of their fundamental groups, namely of the existence of an infinite cyclic normal subgroup. It's now a theorem, of major importance in the understanding of compact 3-manifolds.\\

 We recall here the motivations and applications, the
generalisations and the history of the proof of the Seifert fibre space conjecture.

\section{The Seifert fibre space conjecture and its applications}
\subsection{Reviews on Seifert fibre space} Seifert fibered spaces originally appeared in a paper of Seifert
(1933\footnote{H.\textsc{Seifert}, {\sl Topologie dreidimensionaler gefaserter R\"aume}, Acta Math. 60 (1933),
147-238.\\ \null\qquad An english translation by W.Heil appears in:\\
\null\quad\ \ \ \ H.\textsc{Seifert} and W.\textsc{Threlfall}, {\sl A textbook of topology}, Pure and Applied Math. 89,
Academic Press, 1980.}); they constitute a large of class of 3-manifolds and are totally classified by mean of a finite
set of invariants. They have since widely appeared in the literature for playing a central key-role in the topology of
compact 3-manifolds, and being nowadays (and since the original paper of Seifert) very well known and understood. They
have allowed the developing of promising central concepts in the study of 3-manifolds such as the JSJ-decomposition and
the Thurston's geometrization conjecture.

\subsubsection{Definition of Seifert fibre space}

It's a result of D.Epstein that the Seifert fibered spaces are characterized as those 3-manifolds which admit a
foliation by circles. We rather use as  definition those 3-manifolds which admit a {\sl Seifert
fibration}. In fact our definition is a little more general than the original definition of Seifert,
in order to correctly englobe the case of non-orientable 3-manifolds;
it has now become the modern usual terminology for Seifert fibre spaces.\\

\noindent{\bf Definitions.}\smallskip\\
Let $D^2=\{z\in\C\, |\, |z|\leq 1\}$ stands for the unit disc in $\C$, and $I=[0,1]\subset\R$ for the unit interval.
\smallskip\\
$\bullet$ Let $M$ and $N$ be two 3-manifolds, each being a disjoint union of simple closed curves called {\sl fibres}.
A {\sl fibre-preserving homeomorphism} from $M$ to $N$ is an homeomorphism which sends each fibre of $M$ onto a fibre
of $N$.\smallskip\\
$\bullet$  A \textsl{fibered solid torus} of type $(p,q)$, where $p,q$ are coprime respectively non-negative and
positive integers, is obtained from $D^2\times I$ by identifying $D^2\times 0$ and $D^2\times 1$ by the homeomorphism
which sends $(z,0)$ on $(\mathrm{exp}(2i\pi p/q).z,1)$. Up to fibre-preserving homeomorphism one can (and does) assume
that $0\leq p\leq \frac{q}{2}$; under such conditions solid fibered torus are uniquely determined  by their type
$(p,q)$ up to fibre-preversing homeomorphism (nevertheless all fibered solid torus are clearly
homeomorphic to a solid
torus $D^2\times S^1$).

After the identification $0\times I$ becomes a simple closed curve called {\sl the axis of the solid torus},
 and for
all $z \in D^2\setminus \{ 0\}$,\; $\bigcup_{k=0}^{q-1} \mathrm{exp}(2i\pi kp/q).z\times I$ becomes a simple closed
curve. So that a fibered solid torus is the disjoint union of such simple closed curves, which are called
\textsl{fibres}. The axis of the torus is said to be {\sl exceptional} whenever $p\not=0$ and {\sl regular} otherwise;
in the former case
 $q$ is the {\sl
index} of the exceptional fibre.\smallskip\\
$\bullet$ A {\sl fibered solid Klein bottle} is obtained from $D^2\times I$ by identifying $D^2\times 0$ to $D^2\times
1$ by the (orientation reversing) homeomorphism which sends $(z,0)$ onto $(\bar{z},1)$; it is clearly homeomorphic to
the twisted $I$-bundle over the M\"obius band. A solid Klein bottle is a union of disjoint simple closed curves, called
fibres:
 the image of $z\times I$ for all $z\in\R$
which are called {\sl exceptional fibres} together with the images of $z\times I\cup\bar{z}\times I$ for all
$z\not\in\R$.\smallskip\\
$\bullet$ For a 3-manifold $M$, a {\sl Seifert fibration} is a partition of $M$ into simple closed curves, called
fibres, such that each fibre in the interior $\mathring{M}$ has closed neigbourhood a union of fibres which is
homeomorphic to a solid fibered torus or to a solid Klein bottle $T$ by a fibre-preserving homeomorphism. A fibre of
$M$ which is sent onto an exceptional fibre of $T$ is said to be an {\sl exceptional fibre} (and one can talk of its
index), and is said to be a {\sl regular
fibre} otherwise.\smallskip\\
$\bullet$ A 3-manifold is a {\sl Seifert fibre space} if it admits a Seifert fibration.\smallskip\\
$\bullet$ Let $M$ be a Seifert fibre space; given a Seifert fibration of $M$, if one identifies each fibre to a point,
one obtains a surface $B$, called the {\sl basis} of the fibration. The images of exceptional fibres in $B$ are called
{\sl exceptional points}. The set of exceptional points of $B$ fall into two parts: those coming from the axis of a
solid fibered torus give rise to isolated points in the interior of $B$ (called {\sl conical points}); those coming
from solid fibered Klein bottles are non isolated: they form a closed (non necessarily connected) 1-submanifold in the
boundary of $B$ (its connected components are called {\sl reflector curves}). So that the basis inherits a structure of
2-orbifold
\footnote{I.\textsc{Satake}, {\sl On a generalization of the notion of manifold}, Proc. Nat. Acad. Sciences {\bf 42}
(1956), 359-363.\\ \null\qquad W.\textsc{Thurston}, {\sl The Geometry and Topology of Three-Manifolds}, Princeton
University lecture notes (1978-1981). (In chapter 13.)
}
without corner reflector.\\

The original definition of Seifert, was a little more restrictive, making no considerations of solid fibre Klein
bottles (both definitions agree for orientable 3-manifolds, cause orientability avoids embedded solid Klein bottle).
Modern considerations (see conjectures below, especially Seifert fibre space conjecture) have pointed out that it seems
fairly more natural to enlarge the definition in order to englobe correctly non-orientable 3-manifolds.

In a modern terminology involving the concept of orbifolds, Seifert fibre spaces are those 3-manifolds in the category
of 3-orbifolds which are circle bundles over 2-dimensional orbifolds. One recovers also the original definition of
Seifert by making the assumption that the singular set of the base 2-orbifold only consists of a finite number of cone
points.

\subsubsection{Basic topological properties of Seifert Fibre spaces}

We denote by $S^1$, $S^2$, $S^3$ the spheres of respective dimensions 1, 2 and 3, by $\Bbb K^2$ the Klein bottle, and
by $\Bbb P^2$, $\Bbb P^3$ the
real projective spaces of respective dimensions 2 and 3.\\

With a few well known
exceptions, a Seifert fibered space admits up to fibered preserving homeomorphism a unique Seifert fibration.

\begin{thm} Seifert fibre spaces admit, up to fibre preserving homeomorphism, a unique Seifert fibration, except in the
following cases:\smallskip\\
$\bullet$ a solid torus; it can be fibered with base $D^2$ and one exceptional fibre of arbitrary index.\smallskip\\
(i)\ {\it lens spaces} (including $S^2\times S^1$ and $S^3$); they each have infinitely many fibrations with base $S^2$
and at most 2
exceptional fibres,\smallskip\\
$\bullet$ a $I$-bundle over the Klein bottle; $I\times \Bbb K^2$ has a $S^1$-fibration over the annulus and a Seifert
fibration with base $D^2$ with two exceptional arcs in its boundary; while $I\tilde{\times} \Bbb K^2$ has a fibration
with base $D^2$ and two exceptional fibres
of index 2, and a fibration with base the M\"obius band and no exceptional fibre,\smallskip\\
$\bullet$ the twisted $I$-bundle over the  torus: it is a $S^1$-bundle over the m\"obius band, and has a Seifert
 fibration with base an annulus with one exceptionel boundary component,\\
  (ii)\ prism manifolds; they have a fibration with base $S^2$ and three exceptional fibres of index 2, 2
and $\alpha>1$, and also with base $\Bbb P^2$ and at most one exceptional fibre.\smallskip\\
(iii) the double of the twisted $I$ bundle over the Klein bottle; it has a fibration with base
 $S^2$ and four
exceptional fibres each of index 2, and a fibration with base the Klein bottle and no
exceptional fibre. \smallskip\\
(iv) both $S^2\tilde{\times}S^1$ (the non trivial $S^2$-bundle over $S^1$)  and $\Bbb P^2\times S^1$ have infinitely
many Seifert fibrations, one with base $\Bbb P^2$ and no exceptional fiber, infinitely many with base $\Bbb P^2$ and
one exceptional fiber, one with base a disc with exceptional boundary, and infinitely many with base a disc with
 exceptional boundary and one exceptional fiber.\smallskip\\
 (v) one among three of the four $\Bbb K^2$-bundles over $S^1$; $\Bbb K^2\times S^1$ admits three Seifert fibrations,
 two with base $S^1\times S^1$ or $K^2$ and no exceptional fiber, and one with
 base an annulus with exceptional boundary; $\Bbb KS$ admits three Seifert fibrations,
 two with base
 $K^2$ or $S^1\times S^1$ and no exceptional fiber, and one with base a M\"obius band with exceptional boundary; $S\Bbb
 K$ admits two Seifert fibrations, one with base $\Bbb K^2$ and no exceptional fiber, and one with base a M\"obius band
 with exceptional boundary.\smallskip\\
\indent
 In addition, an homeomorphism between two Seifert fibers spaces which appear
 neither in the above nor in the list below is isotopic to a fibre preserving homeomorphism.\smallskip\\
 $\bullet$ the 3-torus $S^1\times S^1\times S^1$,\smallskip\\
 $\bullet$ the trivial $I$-bundle over the torus $S^1\times S^1\times I$,\smallskip\\
 $\bullet$ the torus-bundle over $S^1$ with characteristic map $-Id$,\smallskip\\
 $\bullet$ the 3-manifold obtained by gluing along the boundary two copies of $I\tilde{\times}\Bbb K^2$
 by an homeomorphism with sends the fibres of one of its two Seifert fibrations onto the fibres of the other.\\

\end{thm}
\noindent {\sl Proof.} The case of closed Seifert fiber spaces has been treated by P.Orlik and
F.Raymond\footnote{P.{\sc Orlik} and F.{\sc Raymond}, {\sl On 3-manifolds with local $SO(2)$ action}, Quart. J. Math.
Oxford, 20 (1969), 143-160.
\smallskip\\ Seifert fiber spaces are those 3-manifolds which
admit a local $SO_2$-action without fixed points.} and one obtains cases (i) to (v). It remains to consider the case of
Seifert fiber spaces with a non-empty boundary. Note that it follows easily from the definition that the boundary only
consists of tori and Klein bottle, which are both a union of fibers.\hfill$\square$

Seifert bundles are clearly foliated by circles; it's a result of D.Epstein \footnote{D.\textsc{Epstein}, {\sl Periodic
Flows on Three-Manifolds}, Annals of Math. 95, 68-82 (1972).} that the converse is true for compact 3-manifolds.
\begin{thm}[D.Epstein]
The Seifert fibre spaces are those compact 3-manifolds which admit a foliation by circles.
\end{thm}

\begin{prop}
The universal cover of a Seifert fiber space has total space homeomorphic to either $\mathbb R^3$, $\mathbb S^2\times
\mathbb R$ or $\mathbb S^3$. Moreover the Seifert fibration lifts in the universal cover to a foliation by lines (in
the former cases) or circles (in the latter case).
\end{prop}

\noindent{\sl Proof.} \hfill$\square$

Recall that a 3-manifold $M$ is {\sl irreducible} if every sphere embedded in $M$ bounds a ball. An irreducible
3-manifold which does not contain any embedded 2-sided $\Bbb P^2$ is said to be $\Bbb P^2$-irreducible. Now Alexander
has shown that $\mathbb S^3$ and $\mathbb R^3$ are irreducible, and it follows that any 3-manifold covered by $\mathbb
S^3$ and $\mathbb R^3$ is $\mathbb P^2$-irreducible. So that the only non $\mathbb P^2$-irreducible Seifert bundles are
covered by $\mathbb S^2\times \mathbb R$, and there are very few such manifolds\footnote{J.\textsc{Tollefson}, {\sl The
compact 3-manifolds covered by $S^2\times R^1$}, Proceedings of the AMS 45 (3), 461-462, (1974)}; one obtains:
\begin{thm}
With the exceptions of $\Bbb P^3\# \Bbb P^3$, $S^2\times S^1$,  $S^2\tilde{\times}S^1$ a Seifert fibre space is $\Bbb
P^2$-irreducible. The only irreducible and non $\Bbb P^2$-irreducible Seifert fibre space is $\Bbb P^2\times S^1$.
\end{thm}

A 3-manifold $M$ is said to be {\sl Haken} if $M$ is $\Bbb P^2$-irreducible and either $M=B^3$ or $M$ contains
an embedded surface $F$, which is :\\
-- {\sl properly embedded}, {\it i.e.} $F\cap\partial M=\partial F$,\\
-- {\sl 2-sided}, {\it i.e.} $F$ has a closed neighborhood in $M$ homeomorphic to $F\times I$,\\
-- {\sl incompressible}, {\it i.e.}, either $F=D^2$ and $\partial F$ is non contractile in $\partial M$ or $F\not=D^2,
S^2$ and the morphism $i_*:\pi_1(F)\longrightarrow \pi_1(M)$ induced by the embedding is injective.

Note that a $\Bbb P^2$-irreducible 3-manifold $M$ whose first homology group $H_1(M,\Z)$ is infinite does contain a
2-sided properly embedded incompressible surface and therefore is Haken; nevertheless there exists Haken 3-manifolds
$\not= B^3$ with finite first homology group (examples can be constructed by performing Dehn obturations on knots).

\begin{thm}
With the exception of lens spaces,  $\Bbb P^3\# \Bbb P^3$, $S^1\times S^2$, $S^1\tilde{\times}S^2$, and $\Bbb P^2\times
S^1$, a Seifert fibre space is either Haken or has base $S^2$ and exactly 3 exceptional fibres; in this last case $M$
is Haken if and only if $H_1(M,\Bbb Z)$ is infinite.
\end{thm}

\begin{thm}
Let $M$ be a 3-manifold equipped with a Seifert fibration.\smallskip\\
\indent (i)\ A regular fibre of $M$ defines up to orientation a homotopy class of loop in $M$; up to orientation
all regular fibres of $M$ define the same homotopy class.\smallskip\\
\indent (ii)\ If the fundamental group $\pi_1(M)$ of $M$ is infinite, a regular fibre of $M$ defines up to conjugacy an
infinite order element of
$\pi_1(M)$.\smallskip\\
\indent In particular, whenever $\pi_1(M)$ is infinite, it contains an infinite cyclic normal subgroup.
\end{thm}

\subsection{The Seifert fibre space conjecture}

\subsubsection{Statement of the Seifert fibre space conjecture}

An infinite fundamental group of a Seifert fibre space contains a normal infinite cyclic subgroup. The Seifert fibre
space conjecture (or SFSC) uses this property to characterize Seifert fibre spaces in the class of oriented irreducible
3-manifolds with infinite $\pi_1$. It can be stated as :

\begin{conj}
Let $M$ be an oriented irreducible 3-manifold whom $\pi_1$ is infinite and contains a non trivial normal cyclic
subgroup. Then $M$ is a Seifert fibre space.
\end{conj}

It generalizes to the non-oriented case by :

\begin{conj}
Let $M$ be a $\Bbb P_2$-irreducible 3-manifold whom $\pi_1$ is infinite and contains a non trivial cyclic normal
subgroup. then $M$ is a Seifert fibre space.
\end{conj}

\noindent{\bf Remarks :} -- Using the sphere theorem together with classical arguments of algebraic topology one shows
that a $\Bbb P_2$-irreducible with infinite $\pi_1$ has a torsion free fundamental group. So that in both conjectures
one can replace :\smallskip\\
"{\sl ... whom $\pi_1$ is infinite and contains a non trivial normal cyclic subgroup.}"\smallskip\\
by :\smallskip\\
"{\sl ... whom $\pi_1$ contains an infinite normal cyclic subgroup.}".\\
-- An oriented Seifert fibre space is either irreducible or homeomorphic to  $S^1\times S^2$ or to $\Bbb P^3\#\Bbb
P^3$. As a consequence of Kneser-Milnor theorem, an oriented non-irreducible 3-manifold $M$ whom $\pi_1$ contains a non
trivial normal cyclic subgroup is either $S^1\times S^2$, or $M'\# C$ with $M'$ irreducible and $C$ simply connected,
or its $\pi_1$ is the infinite dihedral group $\Z_2*\Z_2$. If one accepts the Poincar\'e conjecture (today stated by the
work of Perelman, nominated for a Fields medal in 2006) $M$ is obtained from $S^1\times S^2$, $\Bbb P^3\# \Bbb P^3$ or
from an irreducible 3-manifold by removing a finite number of balls. \\
-- A non-oriented Seifert fibre space is either $\Bbb P^2$-irreducible or $S^1\ltimes S^2$ or $\Bbb P^2\times S^1$. A
large class of irreducible, not $\Bbb P^2$-irreducible 3-manifolds do not admit a Seifert fibration while their $\pi_1$
contains a normal subgroup homeomorphic to $\Z$. We will see that nevertheless the result generalizes in that case by
considering the so-called Seifert fibre spaces mod $\Bbb P$.\medskip\\
\indent The SFSC has now become a theorem, with the (huge) common work of Mathematicians : Waldhausen (\citechrono{wald}),
Gordon and Heil (\citechrono{gh}), Jaco and Shalen (\citechrono{js}), for the Haken case ; Scott (\citechrono{scott}), Mess (\citechrono{mess},
non-published), Tukia (\citechrono{tukia}), Casson and Jungreis (\citechrono{cass}), Gabai (\citechrono{gabai}), for the oriented
non-Haken case ; Heil and Whitten (\citechrono{whit,{heil-whitten}}) for the non-oriented case ; one can also cite Maillot
(\citechrono{maill}, \citechrono{maillot}) and Bowditch (\citechrono{bow2}) who give an alternate proof including the non-published
key-result of Mess.

We can also note that the theorem can be generalized in several ways  : for open 3-manifolds and 3-orbifolds
(\citechrono{maillot}), and for $PD(3)$-groups (\citechrono{bow2}), or by weakening the condition of existence of a normal $\Z$ by
the existence of a non trivial finite conjugacy class (\citechrono{hp}).

\subsection{Motivations}

Three important questions in 3-dimensional topology have motivated the SFSC. First the center conjecture (1960's), then
Scott's strong torus theorem (1978), and finally the Thurston's geometrization conjecture (1980's).
\subsubsection{Center conjecture} It's the problem 3.5 in Kirby's list, attributed to Thurston. :\smallskip\\
\noindent{\bf Conjecture : } {\sl Let $M$ be an oriented irreducible 3-manifold with infinite $\pi_1$ having a non
trivial center, then $M$ is a Seifert fibre space.}\smallskip

Clearly it's an immediate corollary of Seifert fibre space conjecture.
It has first been observed and proved for knot complements (Murasugi (\citechrono{mu}), Neuwirth (\citechrono{ne}) in 1961 for
alternated knots, and Burde, Zieschang (\citechrono{bz}) in 1966 for all knots), and has then been proved in 1967 by
Waldhausen (\citechrono{wald}) in the most general case of Haken manifolds.

An infinite cyclic normal subgroup of the $\P$ of an oriented Seifert fibre space is central if and only if the base of
the fibration is oriented. A Seifert fibre space with infinite $\pi_1$ and with a non-oriented base has a centerless
$\P$ containing a non trivial normal cyclic subgroup. The SFSC generalizes the center conjecture in that sense.

\subsubsection{The torus theorem} The so called "torus theorem" conjecture asserts :\smallskip\\
{\bf Conjecture :} {\sl Let $M$ be an oriented irreducible 3-manifold with $\pi_1(M)\supset \ZZ$. Then either $M$
contains contains an incompressible torus or $M$ is a (small) Seifert.}\smallskip

In has been shown in the case where $M$ is Haken by Waldhausen in 1968 (announced in \citechrono{w2}, written and published
by Feustel in \citechrono{f1,f2}). Note that together with the theory of {\sl sufficiently large} 3-manifolds developed by
Haken during the $60's$ it has finally given rise to the {\sl Jaco-Shalen-Johansen}
decomposition of Haken 3-manifolds (1979).\smallskip\\
\indent
In 1978 Scott proves (\citechrono{stt}) the "strong torus theorem" :\smallskip\\
{\bf Strong torus theorem : } {\sl Let $M$ be an oriented irreducible 3-manifold, with $\pi_1(M)\supset \ZZ$. Then
either $M$ contains an incompressible torus, or $\pi_1(M)$ contains a non trivial normal cyclic subgroup.}\smallskip

With the strong torus theorem, in order to prove the torus theorem it suffices to prove the SFSC.

\subsubsection{The geometrization conjecture}
From the Thurston geometrization conjecture would follow the classification of 3-manifolds.
It conjectures that the pieces obtained in the canonical topological decomposition of an oriented 3-manifold, along
spheres, discs and essential tori, have interiors which admit complete locally homogeneous riemannian metrics. It
follows that their interiors are modelled on one of the 8 homogeneous 3-dimensional geometries : the 3 isotropic ones
(elliptic $\Bbb S^3$, euclidian $\Bbb E^3$, and hyperbolic --the generic-- $\Bbb H^3$), the 2 product ones ($\Bbb
S^2\times\R$ and $\Bbb H^2\times \R$) and the 3 twisted ones ($Nil$, $Sol$ and the universal cover of $SL_2\R$).
Note that an oriented 3-manifold admits a Seifert fibration exactly when its interior is modelled on one the 6
geometries different from $Sol$ and $\H^3$.\\

Thurston ({\sl and al...}) has shown the geometrization conjecture in the Haken case :\smallskip\\
{\bf Thurston geometrization conjecture : } {\sl A Haken 3-manifold $M$
satisfies the Thurston geometrization conjecture.}\\

To prove the conjecture in the remaining cases it suffices to prove the 3 following conjectures for any oriented closed
irreducible 3-manifold $M$ :\smallskip\\
{\sl-1- If $\pi_1(M)$ is finite, then $M$ is elliptic (the othogonalization conjecture).\smallskip\\
-2- If $\pi_1(M)$ is infinite and contains a non trivial normal cyclic subgroup then $M$ is a Seifert fibre space (CSFS.)
\smallskip\\
-3- If $\pi_1(M)$ is infinite and contains no non trivial normal cyclic subgroup then $M$ is geometrizable.} \smallskip\\
Together with the strong torus theorem, the Thurston geometrization theorem and the fact that a closed 3-manifold
modelled on one of the 7 non hyperbolic geometries has a $\pi_1$ containing $\ZZ$,
point 3 becomes :\smallskip\\
{\sl-3- If $\pi_1(M)$ is infinite and contains no $\ZZ$ then $M$ is hyperbolic (hyperbolization conjecture).}
\smallskip\\

So that the SFSC appears as one of the 3 pieces of the geometrization conjecture (usually considered as the easiest
one); it has been the first being proved. The two other pieces have been recently independently proved by the work of
Perelman granted by a Fields medal in 2006 ; Perelman has followed the Hamilton program using the Ricci flow.

\bibliographystylegen{plain}
\bibliographygen{biblio}
\addcontentsline{toc}{subsection}{References}

\section{Historic of the proof}
\subsection{The Haken oriented case} The historic of the proof is the following
:\smallskip\\
{\bf 1967.} Waldhausen (\citechrono{wald}, \citechrono{f1,f2}) shows that a Haken 3-manifold $M$ has a $\pi_1$ with non trivial
center if and only if $M$ is a Seifert fibre space having an oriented base.

It then motivates the SFSC and solves the Haken case when the cyclic normal subgroup is central.\smallskip\\
{\bf 1975.} Gordon and Heil (\citechrono{gh}) show partially SFSC in the Haken case : either $M$ is a Seifert space or it is
obtained by gluing two copies of a non trivial $I$-bundle along a non oriented surface. So that they reduce the
remaining Haken cases to these last 3-manifolds.\smallskip\\
{\bf 1979.} Jaco and Shalen (\citechrono{js}), and independently MacLachlan (non published) achieve the proof for the
remaining Haken oriented 3-manifolds.

\subsection{The non Haken oriented case} Using the Haken case already established, it suffices to restrict
to the closed 3-manifolds.
 the proof has proceeded the following way
:\medskip\\
 \noindent {\bf 1983.} Scott (\citechrono{scott})
(by generalizing a result of Waldhausen in the Haken case) shows that : \\
{\sl \indent Let $M$ and $N$ be 2 closed oriented irreducible 3-manifolds ; where $N$ is a Seifert fibre space with
infinite $\pi_1$. If $\pi_1(M)$ and $\pi_1(N)$ are isomorphic, then $M$ and $N$ are
homeomorphic.}\smallskip\\
\noindent {\bf Remark :} By hypothesis, with the sphere theorem and classic arguments of algebraic topology, $M$ and
$N$ are Eilenberg-MacLane spaces, or $K(\pi,1)$. With the Moise theorem (...), they can be considered in the PL
category. And the condition "$\pi_1(M)$ and $\pi_1(N)$ are isomorphic" can be replaced by "$M$ and $N$ have the same
homotopy type".\smallskip\\
\indent
Hence Scott reduces the Seifert fibre space conjecture to :\smallskip\\
\indent {\sl If $\G$ is the group of a closed oriented irreducible 3-manifold and contains a normal $\Z$, then $\G$ is
the group of a closed oriented Seifert fibre space.}\smallskip

He remarks moreover that $\G$ is the group of a (closed oriented) Seifert fibre space if and only if $\G/\Z$ is the
group of a 2-orbifold (closed eventually non oriented) ; the same result appears in the lemma 15.3 dof \citechrono{bow}. So
finally he reduces the proof of the Seifert fibre space conjecture to the proof of the conjecture :\smallskip\\
\indent{\sl\bf If $\G$ is the group of a closed oriented 3-manifold which contains a normal $\Z$, then $\G/\Z$ is the
group of a (closed) 2-orbifold}.\\

\noindent {\bf Late 80's.} Mess (\citechrono{mess}) shows in a (non published) paper
 :\smallskip\\
\indent{\sl If $M$ is closed oriented and irreducible with $\G=\pi_1(M)$
containing an infinite cyclic subgroup $C$ :\\
-- The covering of $M$ associated to $C$ is homeomorphic to the open manifold $D^2\times S^1$.\\
-- The covering action of $\G/C$ on $D^2\times S^1$ gives rise to an approximatively defined action on $D^2\times 1$ :
following its terminology $\G/C$ is coarse quasi-isometric to the euclidian or hyperbolic plane.\\
-- In the case where $\G/C$ is coarse quasi-isometric to the euclidian plane, then it's the group of a 2-orbifold, and
hence (with the result of Scott) $M$ is a Seifert fibre space.\\
-- In the remaining case where $\G/C$ is coarse quasi-isometric to the hyperbolic plane, $\G/C$ induces an action on
the circle at infinity, which makes $\G/C$ a convergence group.}\smallskip\\
\indent
 A {\sl convergence group} is a group $G$ acting by orientation preserving homeomorphism on the circle,
in such a way that if $T$ denotes the set of ordered triples : $(x,y,z)\in S^1\times S^1\times S^1$, $x\not=y\not=z$,
such that $x,y,z$ appear in that order on $S^1$ in the positive direction, the action induced by $G$ on $T$ is free and
properly discontinuous.\smallskip\\
With the work of Mess, the proof of the SFSC reduces to proving the conjecture :\\
{\sl\bf \indent Convergence groups are groups of 2-orbifolds.}\\

\noindent {\bf 1988.} Tukia shows (\citechrono{tukia}) that some (when $T/\G$ is non compact and $\G$ has no torsion elements
with order $>3$) convergence groups are Fuchsian groups, and in particular
$\pi_1^{orb}$ of 2-orbifolds.\\

\noindent {\bf 1992.} Gabai shows (\citechrono{gabai}) independently from other works and in full generality that convergence
groups are fuchsian groups (acting on $S^1$, up to conjugacy in $Homeo(S^1)$, as restriction on $S^1=\partial \Bbb H^2$
of the natural action of a fuchsian group on $\bar{\Bbb H}^2$).\\

\noindent {\bf 1994.} At the same time Casson and Jungreis show the cases left remained by Tukia (\citechrono{cass}).\\

{\bf Hence the SFSC is proved, and then follow the center conjecture, the torus theorem and one of the 3 geometrization
conjectures.}\\

\noindent{\bf 1999.} Bowditch obtains a different proof of the SFSC (\citechrono{bow}, \citechrono{bow2}). Its proof generalizes to
other kinds of group as  $PD(3)$-groups for example.\\

\noindent{\bf 2000.} In his thesis Maillot, extends the techniques of Mess to establish a proof of the SFSC in the more
general cases of open 3-manifolds and of $3$-orbifolds. This result already known as a consequence of the SFSC (as
proved by Mess and \textsc{al}) and of the Thurston orbifolds theorem has the merit to be proved by using none of
these results and to state as a corollary a complete proof of the SFSC and of the techniques of Mess.

The argument of Mess is pursued in \citechrono{maillot} where he reduces to groups quasi-isometric to a riemanian
plane, and he shows in \citechrono{maill} that they are virtually surface groups (and hence fuchsian groups).

\subsection{The non-oriented case}
The non-oriented case has been treated by Whitten and Heil.\\

\noindent{\bf 1992.} Whitten shows in \citechrono{whit} the SFSC for $M$ non-oriented and irreducible, which is not a fake
$\Bbb P_2\times S^1$, and such that $\pi_1(M)$ contains no $\Z_2*\Z_2$. He obtains in particular the conjecture
conjecture 2, or SFSC in the $\Bbb
P_2$-irreducible case.\\

\noindent{\bf 1994.} Heil and Whitten in \citechrono{heil-whitten} characterise those non oriented irreducible 3-manifolds
which do not contain a fake $\Bbb P_2\times I$ and whom $\pi_1$ contains a normal $\Z$ subgroup and eventually
$\Z_2*\Z_2$ : they are the so called Seifert fibre spaces mod $\Bbb P$ : they are obtained from a Seifert fibered
3-orbifold by **picking of** all neigborhoods of singular points homeomorphic to cones on $\Bbb P_2$ : their boundary
contains in general some $\Bbb P_2$ and otherwise they are merely Seifert fibre spaces. The result becomes :

\begin{thm} Let $M$ be a non oriented irreducible 3-manifold
which does not contain a fake $\Bbb P_2\times I$ with its $\pi_1(M)$ containing a non trivial cyclic subgroup then $M$
is either $\Bbb P_2\times I$ or a Seifert fibre space mod $\Bbb P$.
\end{thm}

They also prove in the non oriented case the torus theorem and deduce the geometrization (i.e. the orientation cover is
geometrizable) whenever the manifold does not contain any fake $\Bbb P_2\times I$.

\subsection{With the Poincar\'e conjecture}\hfill\\

\noindent {\bf 2003-06.} The work of Perelman (2003) states (as a corollary of the orthogonalisation conjecture) the
Poincar\'e conjecture. Hence each simply connected closed 3-manifold is a  3-sphere and there exists no fake ball or fake
$\Bbb P_2\times I$. Hence one can suppress in the SFSC the hypothesis of irreducibility without **alourdir** the
assertion.
\begin{thm} Let $M$ be a 3-manifold whose $\pi_1$ is infinite and contains a non trivial cyclic normal subgroup.
After having filled the spheres in $\partial M$ with balls one obtains either a connected sum of $\Bbb P_2\times I$
with itself or with $\Bbb P_3$, or a Seifert fibre space mod $\Bbb P$.
\end{thm}
One obtains a Seifert fibre space exactly when $\P(M)$ does not contain any $\Z_2$ or in equivalently when $\partial M$
does not contain any $\Bbb P_2$.

\subsection{3-manifold groups with infinite conjugacy classes}\hfill\\

\noindent {\bf 2005.} P. de la Harpe and the author show in \citechrono{hp} that in the case of 3-manifolds and of
$PD(3)$-groups the hypothesis "{\sl contains a non trivial cyclic subgroup}" can be weakened as "{\sl contains a non
trivial finite conjugacy class}" or equivalently as "{\sl has a Von-Neumann algebra which is not a type II-1 factor}".
The proof of these results makes use of the SFSC for 3-manifolds as well as for $PD(3)$-groups (as stated by Bowditch).

\noindent


\addcontentsline{toc}{subsection}{References for the proof}

\noindent
\begin{thebibliography}{mot}\addcontentsline{toc}{subsection}{References for the proof of the Seifert fibre space
conjecture}
%
\bibitem[Se-33]{seifert} H.\textsc{Seifert}, \emph{Topology of 3-dimensional fibered spaces},
Acta Mathematica {\bf 60} (1933), 147--288.\smallskip\\
The english version of the original paper of H.Seifert introducing Seifert fibered spaces (translated from german by W.Heil); among other things he  classifies them and computes their fundamental group.
\bibitem[Mu-61]{mu}
K.\textsc{Murasugi}, \emph{Remarks on torus knots}, Proc. Japan Acad. {\bf 37} (1961), 222.
\bibitem[Ne-61]{ne}
L.\textsc{Neuwirth}, \emph{A note on torus knots and links determined by their groups}, Duke Math. J. {\bf 28} (1961),
545-551.\smallskip\\
Preuve de la conjecture du centre pour les neuds altern\'es.\smallskip
\bibitem[BZ-66]{bz}
G.\textsc{Burde} et H.\textsc{Zieschang}, \emph{Eine kennzeichnung der Torusknotten}, Math. Ann. {\bf 167} (1966),
169-176.\smallskip\\
Preuve de la conjecture du centre pour les compl�ments de noeuds : ceux dont le $\P$ a un centre $\not= 1$ sont
compl�ments de noeuds toriques.\smallskip
\bibitem[Wa-67]{wald}
F.\textsc{Waldhausen}, \emph{Gruppen mit Zentrum und 3-dimensionale Mannigfaltigkeiten}, Topology {\bf 6} (1967),
505-517.\smallskip\\
Les vari�t�s Haken dont le groupe a un centre non trivial sont des fibr�s de Seifert.\smallskip
%
\bibitem[Wa-68]{w2}
F.\textsc{Waldhausen}, \emph{On the determination of some bounded
3-manifolds by their fundamental groups alone},
Proc.Int.Symp.Top.Hercy-Novi, Yugoslavia, (1968), 331-332.\smallskip\\
Annonce du th�or�me du tore dans le cas Haken.\smallskip
%
\bibitem[GH-75]{gh}
C.\textsc{Gordon} et W.\textsc{Heil}, \emph{Cyclic normal subgroups of fundamental groups of 3-manifolds}, Topology
{\bf 14} (1975), 305-309.\smallskip\\
CFS partiellement obtenu dans le cas Haken : soit $M$ est un fibr� de Seifert soit $M$ est obtenu en recollant deux
copies d'un $I$-fibr� non trivial.\smallskip
%
\bibitem[F1-76]{f1}
C.\textsc{Feustel}, \emph{On the Torus theorem and its applications}, Trans. A.M.S. {\bf 217} (1976), 1-43.
%
\bibitem[F2-76]{f2}
C.\textsc{Feustel}, \emph{On the Torus theorem for closed 3-manifolds}, Trans. A.M.S. {\bf 217} (1976),
45-57.\smallskip\linebreak
%
D�monstration du th�or�me du tore dans le cas Haken, comme r�solu par Waldhausen.\smallskip
%
\bibitem[Sc-78]{stt}
P.\textsc{Scott}, \emph{A new proof of the annulus and torus
theorems}, Amer.J.Math. {\bf 2} (1978), 241-277.\smallskip\\
Le "Strong Torus Theorem".\smallskip
%
\bibitem[JS-79]{js}
W.\textsc{Jaco} et P.\textsc{Shalen}, \emph{Seifert fibered spaces in 3-manifolds}, Memoirs A.M.S. {\bf 2}
(1979).\smallskip\\
Fin de  la d�monstration de  CFS dans le cas Haken. Mais aussi surtout le th�or�me Jaco-Shalen-Johansen.\smallskip
%
\bibitem[Sc-83]{scott}
P.\textsc{Scott},
        \emph{There are no fake Seifert fibre spaces with infinite $\pi_1$},
        Annals of Math., {\bf 117} (1983), 35-70.\smallskip\\
Th�or�me de rigidit� pour les fibr�s de Seifert. \smallskip
%
\bibitem[Me]{mess}
G.\textsc{Mess}, \emph{Centers of 3-manifold groups, and groups
which are coarse quasi-isometric to planes}, Unpublished.\smallskip\\
Ram�ne la CFS � prouver que les groupes de convergence du cercle
sont virtuellement des groupes de surface.\smallskip
%
\bibitem[Tu-88]{tukia}
P.\textsc{Tukia}, \emph{Homeomorphic conjugates of Fuchsian
groups}, J.Reine Angew.Math. {\bf 391} (1988), 35-70.\smallskip\\
R�sout partiellement la conjecture des groupes de
convergences.\smallskip
%
\bibitem[Ga-92]{gabai}
D.\textsc{Gabai}, \emph{Convergence groups are Fuchsian groups},
Annals of Math. {\bf 136} (1992), 447-510.\smallskip\\
Montre que les groupes de convergence du cercle sont des groupes
Fuchsiens ; il d�duit en corollaire la CFS, le th�or�me du tore et
la g�om�trisation dans ce cas.\smallskip
%
\bibitem[Wh-92]{whit}
W.\textsc{Whitten}, \emph{Recognizing nonorientable Seifert
bundles}, J.Knot theory \& Ram., {\bf 1} (1992),
471-475\smallskip\\
La CFS dans le cas non-orientable.
%
\bibitem[CJ-94]{cass}
A.\textsc{Casson} et D.\textsc{Jungreis}, \emph{Convergence groups
and Seifert fibered 3-manifolds}, Invent. Math., {\bf 118} (1994),
441-456.\smallskip\\
Apportent ind�pendamment de Gabai une solution en donnant une
solution aux cas laiss�s manquants par Tukia, et en d�duisent la
CFS.\smallskip
%
\bibitem[HW-94]{heil-whitten}
W.\textsc{Heil}, et W.\textsc{Whitten},
        \emph{The Seifert fibre space conjecture and torus theorem for non-orientable 3-manifolds},
        Canad.Math.Bull. {\bf 37} (4) (1994),
        482-489.\smallskip\\
Le cas non-orientable irr�ductible non $\Bbb P_2$-irr�ductible de CFS. Ils en d�duisent le th�or�me du tore et la
g�om�trisation dans ce cas (modulo des faux $\Bbb P_2\times S^1$, c'est � dire modulo la conjecture de
Poincar�.)\smallskip
%
\bibitem[Bo-99]{bow}
B.\textsc{Bowditch}, \emph{Planar groups and the Seifert conjecture}, preprint (1999).\smallskip
\bibitem[Bo-04]{bow2}
B.\textsc{Bowditch}, \emph{Planar groups and the Seifert
conjecture}, J.Reine angew.Math., {\bf 576} (2004), 11-62.\smallskip\\
 Donne une preuve ind�pendante de la CFS ; elle se g�n�ralise aux groupes $PD(3)$.\smallskip
%
\bibitem[Ma-01]{maill}
S.\textsc{Maillot}, \emph{Quasi-isometries of groups, graphs and surfaces}, Comment. Math. Helv., {\bf 76} (1) (2001),
29-60.\smallskip\\
O� il est montr� que les groupes quasi-isom�triques � une surface
riemanienne compl�te simplement connexe sont les groupes
virtuellement de surface.\smallskip
%
\bibitem[Ma-03]{maillot}
S.\textsc{Maillot}, \emph{Open 3-manifolds whose fundamental groups have infinite center, and a torus theorem for
3-orbifolds}, Trans. of the A.M.S., {\bf 355} (11) (2003),
4595-4638.\smallskip\\
Etend la strat�gie de Mess {\it et al}, pour retrouver la CFS dans le cas des 3-orbi�t�s, et en corollaire une preuve
compl�te de CFS.\smallskip
%
\bibitem[HP-05]{hp}
P.\textsc{de la Harpe} et J.-P.\textsc{Pr�aux}, \emph{Alg�bres d'op�rateurs et groupes fondamentaux des 3-vari�t�s},
ArXiV Math.GR/0509449, preprint (2005).\smallskip\\
La condition que le groupe contienne un $\Z$ distingu� dans la CFS peut �tre affaiblie par la condition qu'il contienne
une classe de conjugaison finie non triviale, ou de fa�on �quivalente que son alg�bre de Von-Neumann ne soit pas un
facteur de type $II-1$.

\end{thebibliography}
\end{document}